\documentclass[a4paper, 12pt, twoside, notitlepage]{amsart}

\usepackage{amsmath,amscd}
\usepackage{amssymb}
\usepackage{amsthm}
\usepackage{comment}
\usepackage{graphicx, xcolor}
\usepackage{mathrsfs}
\usepackage[ocgcolorlinks, linkcolor=blue]{hyperref}

\usepackage{bm}
\usepackage{bbm}
\usepackage{url}

\usepackage[utf8]{inputenc}
\usepackage{mathtools,amssymb}
\usepackage{esint}
\usepackage{tikz}
\usepackage{dsfont}
\usepackage{relsize}
\usepackage{url}
\usepackage[shortlabels]{enumitem}
\usepackage{lineno}
\usepackage{enumitem}

\usepackage{verbatim}
\usepackage{dsfont}
\numberwithin{equation}{section}

\allowdisplaybreaks

\graphicspath{{images/}}

\title[Inverse Problems for Twisted Geodesic Flows]{Inverse Problems for Twisted Geodesic Flows}
\author[S.R. Jathar]{Shubham R. Jathar}
\address{ Computational Engineering, School of Engineering Sciences,
Lappeenranta--Lahti University of Technology LUT, Lappeenranta, Finland}
\email {Shubham.Jathar@lut.fi}

\author[J. Railo]{Jesse Railo}
\address{
Computational Engineering, School of Engineering Sciences,
Lappeenranta--Lahti University of Technology LUT, Lappeenranta, Finland}
\email{Jesse.Railo@lut.fi}

\begin{document}

\begin{abstract}
The article surveys inverse problems related to the twisted geodesic flows on Riemannian manifolds with boundary, focusing on the generalized ray transforms, tensor tomography, and rigidity problems. The twisted geodesic flow generalizes the usual geodesic flow, allowing to model richer dynamics.
\medskip

\noindent{\bf Keywords.} generalized geodesic flows; integral geometry; inverse problems.

\noindent{\bf Mathematics Subject Classification (2020)}: Primary 44A12, Secondary 58C99; 37E35.
\end{abstract}
	\maketitle

\section{Twisted Geodesics}

Let \((M,g)\) be a Riemannian manifold of dimension \(n \geq 2\) with smooth boundary, and let \(\nabla\) denote the Levi--Civita connection. 
Consider the family of twisted geodesics \(\gamma\) for \((x,\xi) \in TM\), where each \(\gamma\) satisfies the geodesic equation
\[
\nabla_{t} \dot{\gamma}(t) = G(\gamma(t),\dot{\gamma}(t))-\frac{\langle G(\gamma(t), \dot{\gamma}(t)),\dot\gamma(t) \rangle}{|\dot \gamma(t)|^2}\dot \gamma(t), \,\, (\gamma(0),\dot{\gamma}(0)) = (x,\xi),
\]
where \(\nabla_t\) denotes the covariant derivative along \(\dot{\gamma}\), and \(G \in C^{\infty}(TM)\). Depending on the nature of the vector field \(G\), special types of twisted geodesics arise:

\begin{itemize}
    \item \textbf{Geodesics:} When \(G = 0\), the equation describes a standard geodesic, representing the path of a freely falling particle in the absence of external forces.
    
    \item \textbf{Magnetic curves:} When \(G(\gamma,\dot{\gamma}) = Y_{\gamma}(\dot{\gamma})\), where \(Y\) is the Lorentz force associated with a closed 2-form \(\Omega\) on \(M\), modeling a magnetic field, and \(Y: TM \rightarrow TM\) is the unique bundle map satisfying \(\Omega_x(\xi,\eta) = \langle Y_x(\xi), \eta \rangle_g\) for all \(x \in M\) and \(\xi, \eta \in T_xM\).
    
    \item \textbf{Gaussian thermostats:} When \(G(\gamma,\dot{\gamma}) = E(\gamma) - \frac{\langle E(\gamma), \dot{\gamma} \rangle_{\gamma}}{|\dot{\gamma}|^2}\dot{\gamma}\), where \(E\) is a vector field on \(M\). This concept was introduced by Hoover as a tool for studying dynamical systems in statistical mechanics.
    
    \item \textbf{\(\mathcal{MP}\) -geodesics:} When \(G(\gamma,\dot{\gamma}) = Y_{\gamma}(\dot{\gamma}) - \nabla U(\gamma)\), where \(U\) is a smooth function on \(M\) and \(Y\) is the Lorentz force as in the magnetic case. This type appears in the context of mechanics.
\end{itemize}
The twisted geodesic flow preserves unit speed, i.e., if \(|\dot{\gamma}(0)| = 1\), then \(|\dot{\gamma}(t)| = 1\) for all \(t\). Let \(SM\) denote the unit sphere bundle of \(M\). For \(t \geq 0\) and \((x,v) \in SM\), the twisted geodesic flow \(\phi_t: SM \to SM\) is defined as \(\phi_t(x,\xi) = (\gamma(t), \dot{\gamma}(t))\), where \(\gamma(t)\) is the unique twisted geodesic with initial conditions \(\gamma(0) = x\) and \(\dot{\gamma}(0) = \xi\). 
Define \(\tau(x,\xi)\) as the travel time, the first time \(\gamma_{x,\xi}\) hits \(\partial M\). A manifold is called nontrapping if \(\tau(x,\xi) < \infty\) for any \((x,\xi) \in SM\).

For recent developments in the geodesic case, we refer to \cite{Ilmavirta:Monard:2019} and the textbooks \cite{GIP2D,Sharafutdinov}. This survey limits to consider the cases of twisted geodesics only. We regret that related questions on closed manifolds are out of scope (for example, see \cite{assylbekov2021att}).
\section{Tensor Tomography}

A central question in integral geometry is whether a function can be determined by its integrals over a family of curves. The twisted ray transform \(I f: SM \rightarrow \mathbb{R}\) of a real-valued function \(f \in C(SM)\) is defined by
\[
I f(x, v) := \int_0^{\tau(x, v)} f\left(\phi_t(x, v)\right) \, \mathrm{d}t.
\]
The twisted ray transform acting on symmetric \(m\)-tensor fields is given by
\[
I_m: C^{\infty}\left(S^m\left(T^* M\right)\right) \rightarrow C^{\infty}\left(\partial_{+} SM\right),
\]
\[
I_m f(x,\xi) = \int_0^{\tau(x,\xi)} f_{i_1 \ldots i_m}(\gamma_{x,\xi}(t)) \dot\gamma_{x,\xi}^{i_1}(t) \ldots \dot\gamma_{x,\xi}^{i_m}(t) \, \mathrm{d}t,
\]
for \((x, \xi) \in \partial_{+} SM\).
The kernel of this transform is nontrivial.

The generator of the magnetic flow \(\textbf{G}_{\mu}\) is given by \cite{Dairbekov-Paternain-SU-2007-magnetic-rigidity}
\[
\textbf{G}_{\mu}(x,\xi) = X(x,\xi) + Y_i^j(x) \xi^i \frac{\partial}{\partial \xi^j},
\]
where \(X\) is the generator of the geodesic flow. Similarly, the generator of the \(\mathcal{MP}\)-flow is \cite{mthon2024lbr} given by
\[
\textbf{G}_{\mathcal{MP}}(x, \xi) = \textbf{G}_\mu(x, \xi) - g^{ij} \partial_{x^i} U(x) \frac{\partial}{\partial \xi^j}.
\]
For surfaces, the generator of a twisted flow is given by \cite[Lemma 7.4]{Marry:Paternain:2011:notes}
\[
\textbf{G}(x,\xi) = X(x,\xi) + G(x,\xi)V(x,\xi),
\]
where \(V\) is a vertical vector field.

An element \(f\) is called a potential if \(f = \textbf{d} w\), where the operator \(\textbf{d}\) depends on the flow. A pair \([h, \beta] \in \mathbf{L}^2(M)\) is called a magnetic potential if there exists \([v, \varphi] \in \mathcal{H}_0^1(M)\) such that
\[
\binom{h}{\beta} = \mathbf{d} \binom{v}{\varphi}, \quad \mathbf{d} := \left(\begin{array}{cc} d^s & 0 \\ -Y & d \end{array}\right),
\]
where \(d^s\) is the symmetric derivative on covector fields, and \(d\) is the usual differential on functions \cite[Def. 3.5]{Dairbekov-Paternain-SU-2007-magnetic-rigidity}. In the case of the \(\mathcal{MP}\)-geodesic flow, the operator \(\textbf{d}\) is given by
\[
\textbf{d} = \begin{pmatrix} d^s & 0 & -\frac{1}{2} g \\ -Y & d & 0 \\ -(d U, \cdot)_g & 0 & k-U \end{pmatrix},
\]
and \(w = [v, \varphi, \eta]\) \cite[Definition 4.3]{mthon2024lbr}. For general twisted geodesics, the potential takes the form \(\textbf{G}w\) for some \(w \in C^{\infty}(SM)\) where $\textbf{G}$ is the generator of the flow $\phi$.

The twisted ray transform \(I_m\) is called solenoidally injective ({s-injec}-tive) for \(m \geq 0\) if its kernel consists only of the corresponding potentials of degree $m-1$ and vanishing on the boundary $\partial SM$. We note that this concept makes mainly sense when $m=0,1$ or $\textbf{G}$ maps the spherical harmonics of degree $m$ to the degree $m+1$ and lower. Below are some known injectivity results for the twisted ray transform and its special cases:

\begin{enumerate}[label=(\roman*)]

    \item For simple twisted flows on Finsler surfaces, the injectivity of the twisted ray transform $I_m$ for $m=0,1$ is proved in \cite[Theorem 1.1]{Assylbekov:Dairbekov:2018}. These results were generalized to certain broken ray transforms under negative curvature assumptions in \cite{RefJ21}.
     \item  For \(m \geq 2\), s-injectivity is shown on simple thermostatic surfaces with non-positive thermostat curvature \cite{Assylbekov:Zhou:2017}.
    \item In \cite[Theorem 2.2]{MR2365669}, the injectivity of the ray transform over an analytic regular family of curves is proved using analytic microlocal arguments, together with a stability result for the localized normal operator \cite[Theorem 2.3]{MR2365669}.
    
    \item Using Melrose’s scattering calculus, the local injectivity of the ray transform for general curves on manifolds of dimension \(d \geq 3\) is proved in the Appendix of \cite{Uhlmann:Vasy:2016} for \(m = 0\). A simlar result for the local magnetic ray transform with \(m = 1, 2\) was proved later in \cite{Zhou:2018:local}.
    
    \item It is demonstrated that the local problem for the twisted ray transform on a surface with conjugate points is ill-posed \cite{MR4572203} by considering the recovery of singularities.
    
    \item For a compact, nontrapping Riemannian surface with a strictly convex boundary with respect to the twisted flow and real analytic $(M,g,SM,\partial_+SM, \textbf{G})$, the twisted ray transform \(I_0\) is injective \cite[Theorem 1.4]{mazzucchelli2023general}.
\end{enumerate}

\section{Matrix Attenuated and Nonabelian Twisted Ray Transforms}

For \((x,v) \in SM\), consider the map \(U_{x,v}:[0,\tau(x,v)]\to \mathbb{C}^{n\times n}\) that satisfies the matrix ordinary differential equation (ODE) along a twisted curve \(\gamma_{x,v}\):
\[
\dot{U}_{x,v}(t) + \mathcal{A}(\gamma_{x,v}(t), \dot{\gamma}_{x,v}(t))U_{x,v}(t) = 0, \quad 
U_{x,v}(\tau(x,v)) = \operatorname{Id},
\]
where \(\mathcal{A} \in C^{\infty}(SM; \mathbb{C}^{n\times n})\). The nonabelian twisted ray transform (or twisted scattering data) of \(\mathcal{A}\) is defined as \(C_{\mathcal{A}}:\partial_+SM\to GL(n,\mathbb{C})\), where \(C_{\mathcal{A}}(x,v) = U_{x,v}(0)\). Additionally, consider the unique solution \(u(t)\) to the vector-valued ODE
\[
\dot{u} + \mathcal{A}(\gamma_{x,v}(t), \dot{\gamma}_{x,v}(t)) u = -f\left(\gamma_{x, v}(t), \dot{\gamma}_{x, v}(t)\right), \quad u(\tau(x, v)) = 0,
\]
where \(\gamma_{x, v}(t)\) is a twisted geodesic. The attenuated twisted \(X\)-ray transform of \(f\) is then defined as
\[
I_{\mathcal{A}}^\lambda(f)(x, v) := u(0),
\]
for \((x, v) \in \partial_{+} SM\). There is a relationship between the nonabelian and attenuated twisted ray transforms, expressed by the pseudolinearization identity
\[
I_{E(\mathcal{A}, \mathcal{B})}^\lambda(\mathcal{A}-\mathcal{B}) = C_{\mathcal{A}}^\lambda\left[C_{\mathcal{B}}^\lambda\right]^{-1} - \operatorname{Id}, \quad E(\mathcal{A}, \mathcal{B}) U:=\mathcal{A} U - U \mathcal{B},
\]
as stated for surfaces in \cite[Lemma 3.7]{jathar2023loop}. This identity extends to higher-dimensional manifolds with analogous proofs. If the attenuated transform is injective, this implies injectivity of the nonabelian transform via the pseudolinearization identity. We refer to the survey \cite{MR4059257} for more details of the nonabelian ray transform and its applications.

Let us consider a special form of attenuation \(\mathcal A(x,v) = A_x(v) + \Phi(x)\), where \(A\) is an \(n\times n\) matrix, each entry being a smooth one-form on \(M\) (a connection), and \(\Phi\) is an \(n\times n\) matrix, each entry being a smooth function on \(M\) (a Higgs field). The natural kernel of the attenuated twisted ray transform includes elements of the form \(\left(\mathbf{G}+A+\Phi\right) p\), where \(p \in C^{\infty}\left(SM; \mathbb{C}^n\right)\) and \(p|_{\partial(SM)}=0\). Below are known results related to these transforms:

\begin{enumerate}[label=(\roman*)]
    \item If \((M,g)\) is a simple magnetic surface, for skew-Hermitian attenuation pairs \((A, \Phi)\), the attenuated and nonabelian magnetic ray transform is injective up to the gauge \cite{Ainsworth:2015}. This is generalized to arbitrary attenuation pairs \((A, \Phi)\) in \cite{jathar2023loop}.
    
    \item If \((M,g)\) is a simple thermostatic surface with negative thermostatic curvature, the attenuated and nonabelian thermostatic ray transform for arbitrary attenuation pairs \((A, \Phi)\) is injective up to the gauge \cite{assylbekov2021att}.
\end{enumerate}

\section{Boundary, Scattering, and Lens Rigidity}

The problems of boundary, scattering, and lens rigidity have primarily been studied in the context of geodesic, magnetic, and \(\mathcal{MP}\)-systems. It is known that magnetic and \(\mathcal{MP}\)-geodesics minimize the corresponding time-free action \cite{Assylbekov-Zhou-2015-magnetic-rigidity-potential,Dairbekov-Paternain-SU-2007-magnetic-rigidity}. Using this time-free action, one can define the boundary action function. The scattering relation \(\mathcal{S}: \partial_{+} SM \rightarrow \partial_{-} SM\) is defined by \(\mathcal{S}(x,\xi) = (\gamma_{x,\xi}(\tau(x,\xi)), \dot{\gamma}_{x,\xi}(\tau(x,\xi)))\) for \((x,\xi)\in \partial_+(SM)\). The lens data consists of the pair \((\mathcal{S}, \tau|_{\partial_{+} SM})\). The following problems have been of great interest:
\begin{itemize}
    \item \textbf{Boundary Rigidity Problem:} Can the metric and \(\Omega\) (or \(\Omega\) and \(U\) for \(\mathcal{MP}\)-systems) be determined, up to a gauge, from the boundary action function?
    \item \textbf{Scattering Rigidity Problem:} Can the metric and \(\Omega\) (or \(\Omega\) and \(U\) for \(\mathcal{MP}\)-systems) be determined, up to a gauge, from the scattering data?
    \item \textbf{Lens Rigidity Problem:} Can the metric and \(\Omega\) (or \(\Omega\) and \(U\) for \(\mathcal{MP}\)-systems) be determined, up to a gauge, from the lens data?
\end{itemize}

The following advances related to these rigidity problems have been obtained:
\begin{enumerate}[label=(\roman*)]
    \item Boundary and scattering rigidity have been established for several classes of metrics, including simple conformal metrics, simple real-analytic metrics, and all simple magnetic surfaces\\\cite{Dairbekov-Paternain-SU-2007-magnetic-rigidity}. Similar results for \(\mathcal{MP}\)-systems are proven in \cite{Assylbekov-Zhou-2015-magnetic-rigidity-potential} using boundary action functions at two energy levels, and in \cite{mthon2023bsr} for a single energy level with a new notion of gauge. Magnetic scattering rigidity for real-analytic manifolds is shown in \cite{Herros-Vargo-2011}.
    
    \item Reconstruction of \(g\) (where \(g=cg_0\) is a scalar function times a known metric) and \(\Omega\) from the scattering relation for simple surfaces is demonstrated in \cite{MR2671103}.
    
    \item As an application of magnetic boundary rigidity, it is shown that surfaces of constant curvature cannot be altered in a small region while maintaining all curves of fixed constant geodesic curvature closed \cite{Herreros2012-magnetic-scattering-rigidity}.
    
    \item The partial data magnetic lens rigidity problem for conformal metrics is solved using Melrose's scattering calculus on manifolds of dimension \(n \geq 3\) in \cite{MR3917833}.
\end{enumerate}


\section*{Acknowledgements}
We thank Sean Holman and Bill Lionheart, the organizers of the minisymposium M13: Applications of Rich Tomography at the IPMS 2024, Malta, for giving us the opportunity to present our work. The work of S.R.J. and J.R. was supported by the Research Council of Finland through the Flagship of Advanced Mathematics for Sensing, Imaging and Modelling (decision number 359183).

\bibliography{math}

\bibliographystyle{alpha}
\end{document}